\documentclass{amsproc}

\usepackage{a4wide}
\usepackage{amsmath}
\usepackage{enumerate}

\usepackage{amsmath,amsthm}
\usepackage{amsfonts}
\usepackage{amssymb}

\usepackage[T1]{fontenc}

\newtheorem{rem}{Remark}{\bf}{\rm}{\rm}

\newtheorem*{theorem*}{Theorem}
\newtheorem*{cor*}{Corollary}

\newtheorem{theorem}{Theorem}
\newtheorem{cor}{Corollary}

\newtheorem{lem}{Lemma}

\newtheorem{ex}{Example}{\rm}{\rm}

\def\Real{\mathbb{R}}

\def\Vol{\text{\rm Vol}}

\newcommand{\mathR}{\mathbb{R}}
\newcommand{\frakso}{\mathfrak{so}}

\newcommand{\g}{\mathfrak{g}}

\newcommand{\so}{\mathfrak{so}}

\newcommand{\zr}{\ltimes}

\newcommand{\Su}{\mathfrak{S}}

\def\id{\mathop\text{\rm id}\nolimits}

\usepackage{stackrel}

\newcommand{\be}{\begin{equation}}
\newcommand{\ee}{\end{equation}}

\let\leq=\leqslant
\let\geq=\geqslant

\begin{document}

\title{Lorentzian connections with parallel twistor-free torsion}

\author{Igor Ernst}\thanks{$^1$Department of Mathematics and Statistics, Masaryk University, Faculty of Science, Kotl\'a\v{r}sk\'a 2, 611 37 Brno, Czech Republic}

\author{Anton S. Galaev}\thanks{$^2$University of Hradec Kr\'alov\'e, Faculty of Science, Rokitansk\'eho 62, 500~03 Hradec Kr\'alov\'e,  Czech Republic\\
E-mail: anton.galaev(at)uhk.cz}

\begin{abstract}

We describe Lorentzian manifolds that admit metric connections with parallel torsion having zero twistorial component and non-zero vectorial component. We also describe Lorentzian manifolds admitting metric connections with closed parallel skew-symmetric torsion.

\vskip0.1cm

\vskip0.1cm



\end{abstract}

\maketitle


\section{Introduction}

Metric connection with torsion appear both in the context of differential geometry and mathematical physics. Under the action of the pseudo-orthogonal group, the torsion of a metric connection on a pseudo-Riemannian manifold  decomposes into three components: vectorial, twistorial, and skew-symmetric. Riemannian connections with skew-symmetric torsion have been studied in many works,  e.g.,  \cite{AF04,A10,CMS21}.
Connections with skew-symmetric torsion appear in certain supergravity theories, see, e.g., \cite{FFSPPM,FI02,MSh,Str}. We studied Lorentzian metric connections with parallel skew-symmetric torsion in
\cite{EG22}. 

Recently Moroianu and Pilca \cite{MP} classified complete simply connected Riemannian manifolds that admit metric connections with parallel torsion having zero twistorial component and non-zero vectorial component. It is shown that each such manifold is a warped product of the line and a complete Riemannian manifold admitting a parallel 3-form $\tau$ satisfying the condition $\tau(X)\cdot\tau=0$ for all vector fields $X$, see Theorem \ref{ThMP5.3} below.
Then complete Riemannian manifolds carrying such  3-forms were classified, see Theorem \ref{ThMP5.2} below.

\emph{In the present paper we consider the corresponding problem in the Lorentzian signature}. The first result (Theorem \ref{ThComplg}) states that if a Lorentzian manifold $(M,g)$ admits a metric connection with parallel torsion having zero twistorial component and non-zero vectorial component, then $(M,g)$ cannot be geodesically complete, and even more it cannot  satisfy no weaker condition of geodesic completeness:   space-like,  time-like,  or light-like geodesic completeness. 
We distinguish  two cases depending on  the vector field $\xi$ defining the vectorial component of the torsion: $\xi$ is isotropic and $\xi$ is not isotropic. If $\xi$ is non-isotropic and complete, then Theorem \ref{Th_nonis} provides a result similar to Theorem \ref{ThMP5.3}. Next, Theorem \ref{ThLordT=0} describes Lorentzian manifolds admitting parallel 3-forms $\tau$ satisfying the above condition. \emph{In fact, Theorem \ref{ThLordT=0} describes Lorentzian manifolds admitting metric connections $\nabla$ with closed $\nabla$-parallel skew-symmetric torsion}.
 Next, we assume that the vector field $\xi$ is isotropic. If the norm of the skew-symmetric component $S$ of the torsion is zero, then the manifold under consideration may be described as a manifold $(M_0,g_0)$ admitting a parallel isotropic vector field and a parallel 2-form on the corresponding screen bundle (Theorem \ref{Th_isotropic_deg}).
Suppose that the norm of $S$ is non-zero. If the dimension of the manifold is 3, then we prove in  Theorem \ref{Th_isotropic_nondeg_dim=3}
that the manifold is a  Kundt space of a special form. If the dimension is greater than~3, then the manifold under consideration is a locally warped product of the manifold $(M_0,g_0)$ and a Riemannian manifold again admitting a 3-form with the above properties.

A large class of metric connections with torsion appear as homogeneous structures.
A homogeneous structure on a pseudo-Riemannian manifold is a metric connection $\nabla$ with $\nabla$-parallel curvature and torsion. By the Ambrose-Singer Theorem, a complete pseudo-Riemannian manifold is reductive homogeneous if and only if it admits a homogeneous structure \cite{CLBook,TV}.  Let us mention several results about homogeneous structures related to the results of the present paper. 
In~\cite{GO97} it is shown that Lorentzian space forms do not admit homogeneous structures of vectorial type. Next, a homogeneous Lorentzian space admitting a homogeneous structure of isotropic vectorial type is  a singular homogeneous plane wave~\cite{MA01,BOL}.  The paper \cite{Meessen} provides a classification of Lorentzian homogeneous structures with the 
torsion having zero twistorial component and non-zero vectorial component.
These structures are exhausted by locally symmetric spaces and singular homogeneous plane waves. We thus generalize these results to general connections with parallel twistor-free torsion having non-zero vectorial component.

\section{Preliminaries}

Let $(M,g)$ be a pseudo-Riemannian manifold. Denote by $\nabla^g$ the Levi-Civita connection on $(M,g)$. A connection $\nabla$ on $(M,g)$ is called metric if $\nabla g=0$.  Using the metric $g$ we identify the tangent and cotangent bundles, and we use the obvious identifications for tensors. In particular, 
we identify a field of bivectors $X \wedge Y$ with the field of endomorphism 
$$ (X \wedge Y) Z = g(X,Z)Y - g(Y,Z)X,$$ and for a 3-tensor $B$ by abuse of notation we will write
$$ B(X, Y, Z) = g(B(X, Y), Z)=g(B(X)Y,Z).$$

Let $A\in \Gamma(TM)\otimes (\wedge^2 \Gamma(TM))$. Then the equality 
$$\nabla_XY=\nabla^g_XY+A(X,Y),\quad X,Y\in \Gamma(TM),$$
defines the metric connection $\nabla$ with the torsion
$$T(X,Y,Z)=A(X,Y,Z)-A(Y,X,Z),$$ see, e.g, \cite{A10}. Let $\xi\in \Gamma(TM)$ and  $S\in \Omega^3(M)$, then the tensor  
$$A(X)=X \wedge \xi + \frac{1}{2} S(X)$$ defines the metric connection 
\begin{equation} \label{connectionnabla}
\nabla_X = \nabla_X^g + X \wedge \xi + \frac{1}{2} S(X), 
\end{equation}
with the torsion 
$$ T(X,Y) = (X \wedge Y)\xi + S(X,Y).$$
The torsion $T$ has vectorial component $\xi$, skew-symmetric component $S$ and zero twistorial component. In particular, a 3-form $\tau$ defines the metric connection 
\begin{equation} \label{connectionnablaTAU}\nabla_X = \nabla_X^g  + \frac{1}{2} \tau(X)\end{equation}
with skew-symmetric torsion $\tau$.

Let $\nabla$ be given by \eqref{connectionnabla}. Suppose that  $\nabla T=0$. It is clear that this is equivalent to the conditions $\nabla \xi=0$ and $\nabla S=0$.
Consider the first Bianchi identity
\[ \underset{X Y Z}{\Su}R(X,Y)Z = \underset{X Y Z}{\Su}\{T(T(X,Y),Z) + (\nabla_X T)(Y,Z) \}, \] where $\underset{X Y Z}{\Su}$ denotes the cyclic sum with respect to $X, Y, Z$. It is easy to check that the identity may be rewritten in the form
\begin{equation}\label{Bidcomp}
\underset{X Y Z}{\Su}R(X,Y)Z = \underset{X Y Z}{\Su}S(S(X,Y),Z) + \underset{X Y Z}{\Su}g(X,\xi)S(Y,Z) + \underset{X Y Z}{\Su}S(\xi,Y,Z)X.
\end{equation}

For a 3-form $\tau$ define the 4-form $\sigma_\tau$ by the equality
$$\sigma_\tau(X,Y,Z)=\underset{X Y Z}{\Su}\tau(\tau(X,Y),Z).$$
 It holds that
\begin{equation}\label{T(X)T}
(\tau(X)\cdot \tau)(Y,Z,V)=-\sigma_\tau(X,Y,Z,V).
\end{equation}
Consider the connection \eqref{connectionnablaTAU}. It holds
\begin{equation}\label{dT1}
d\tau(X,Y,Z,V) =\underset{X Y Z}{\Su}\{(\nabla_X\tau)(Y,Z,V)\}-(\nabla_V\tau)(X,Y,Z)+2\sigma_\tau(X,Y,Z,V), 
\end{equation}
see, e.g., \cite{FI02}. 
From \eqref{connectionnablaTAU} it follows that
\begin{equation}\label{dT2}
\nabla\tau=\nabla^g \tau-\frac{1}{2}\sigma_\tau.
\end{equation}
If $\nabla\tau=0$, then
 the first Bianchi identity for the connection $\nabla$  may be written in the form
\begin{equation}\label{Bident}\underset{X Y Z}{\Su}R(X,Y)Z=\sigma_\tau(X,Y,Z).\end{equation}

Note that by \eqref{T(X)T}, the condition $\sigma_\tau=0$ is equivalent to the condition $\tau(X)\cdot \tau = 0$ for all $X \in \Gamma(TM)$.
From  \eqref{dT1} and \eqref{dT2} easily follows 

\begin{lem}\label{lem} Let $(M,g)$ be a pseudo-Riemannian manifold and $\nabla$ a metric connection with skew-symmetric torsion $\tau$ on it. Then the following conditions are equivalent 
	\begin{itemize}
		\item[1.] $\nabla \tau=0$ and $d\tau=0$;
		\item[2.] $\nabla \tau=0$ and $\sigma_\tau= 0$;
		\item[3.] $\nabla^g \tau=0$ and $\sigma_\tau= 0$.
	\end{itemize}
\end{lem}

The following theorem from \cite{MP} describes complete simply connected  Riemannian manifolds with parallel twistor-free torsion.

\begin{theorem}\label{ThMP5.3}
	A complete simply connected Riemannian manifold $(M, g_M)$ carries a metric connection with parallel twistor-free torsion if and only if $(M, g_M)$ is homothetic to a warped product $( \mathbb{R}\times N,  (dt)^2+e^{2t}g_N)$, where $(N, g_N)$ is a complete simply connected Riemannian manifold carrying a parallel 3-form $\tau \in \Omega^3(N)$ such that $\sigma_\tau=0$.
\end{theorem}

In \cite{AFeF15} it is shown that an irreducible complete simply connected Riemannian manifold of dimension greater than $4$ with parallel skew torsion $\tau$ satisfying $\sigma_\tau=0$  is a simple compact Lie group with bi-invariant metric or its dual noncompact symmetric space.
The following is a more general theorem from~\cite{MP}.

\begin{theorem}\label{ThMP5.2}
	Let $(N, g)$ be a complete simply connected Riemannian manifold carrying a metric connection with parallel skew-symmetric torsion $\tau$ which satisfies $\sigma_\tau=0$. Then $(N,g)$ is one of the following:
	\begin{itemize}
		\item[1.] $(N,g)$ is an oriented 3-dimensional Riemannian manifold and $\tau$ is a constant multiple of its Riemannian volume form;
		\item[2.] $(N,g)$ is a simple Lie algebra with an $\operatorname{ad}$-invariant metric $g$ and $\tau$ is a constant multiple of its canonical 3-form;
		\item[3.] $(N,g)$ is an irreducible symmetric space of type II or type IV and $\tau$ is a constant multiple of its canonical 3-form;
		\item[4.] $(N,g)$ is a Riemannian manifold and $\tau = 0$;
		\item[5.] $(N,g)$ is a product of some of the manifolds as in 1--4.
			\end{itemize}
\end{theorem}

Now we recall the definition of Kundt and Walker spaces.

A Kundt space is a Lorentzian manifold $(M,g)$ admitting an isotropic vector field $p$ satisfying the following conditions:
\begin{equation}\label{Kundtfield}
\nabla^g_pp=0,\quad \operatorname{tr}\nabla^gp=0,\quad ||(\nabla^gp)^{sym}||_g=0, \quad ||(\nabla^gp)^{alt}||_g=0,
\end{equation}
i.e., the vector field $p$ is  geodesic, expansion-free, shear-free and twist-free,
see, e.g., \cite{CHPP}. Locally there exist coordinates $v,x^1,\dots x^n,u$ such that $g$  takes the form 
\begin{equation} \label{Kundt}
g = 2 dv du + h + 2A du + H (du)^2, 
\end{equation}
where $h=\sum_{i,j=1}^nh_{ij}(x^1, \dots, x^n, u) dx^i dx^j$ is a $u$-family of local Riemannian metrics,\\ $A=\sum_{i = 1}^n A_i(v,x^1, \dots, x^n, u) dx^i$ is a 1-form, and $H=H(v,x^1, \dots, x^n, u)$ is a local function. The vector field $p$ is proportional to $\partial_v$.

A Walker manifold (see, e.g.,  \cite{Walkerbook}) is a Lorentzian manifold $(M,g)$ admitting a $\nabla^g$-parallel distribution of isotropic lines $\ell$. Locally $g$ is given by \eqref{Kundt} with the functions $A_i$ independent of $v$. If moreover the induced connection in $\ell$ is flat, i.e.,   $(M,g)$ admits (local) $\nabla^g$-parallel isotropic vector field $p$, then in the function $H$ may be chosen to be  independent of $v$.

Let $\nabla$ be a metric connection on a Lorentzian manifold $(M,g)$ and suppose that $p$ is a $\nabla$-parallel isotropic vector field on $M$.
Since the metric $g$ is $\nabla$-parallel, the distribution $p^\bot$ is $\nabla$-parallel.
The bundle $E=p^\bot/\left<p\right>$ is called  \emph{the screen bundle}, see \cite{LScr,EG22} for  details.    There is the obvious projection $p^\bot\to E$. The connection $\nabla$ induces a connection on $E$ (which we denote by the same symbol): if $X$ is a vector field on $M$ and $Y$ is a section of $E$, then  $\nabla_XY$ is the projection to $E$ of the vector field $\nabla_X \tilde Y$, where  $\tilde Y$ is an arbitrary section of $p^\bot$ such that its projection to $E$ is $Y$.

We will consider a Witt basis $p,e_1,\dots, e_n,q$ of the Minkowski space $\Real^{1,n+1}$. The non-zero values of the Minkowski metric with respect to such basis are $g(p,q)=g(q,p)=g(e_i,e_i)=1$. We will denote by $\Real^n$ the Euclidean subspace of $\Real^{1,n+1}$ spanned by the vectors
$e_1,\dots, e_n$. With respect to the basis $p,e_1,\dots,e_n,q$, the subalgebra of $\so(1,n+1)$ preserving the isotropic line $\Real p$ has the following matrix form:
$$
\so (1, n + 1)_{\Real p} = \left\{
\left. 
\begin{pmatrix}     
a & -X^t & 0 \\
0 & A & X \\
0 & 0 & -a
\end{pmatrix} \right|
\begin{matrix}     
a \in \mathR \\
A \in \frakso (n) \\
X \in \mathR^n
\end{matrix} \right\} .
$$	
The above matrix is identified with the bivector
$$-ap\wedge q+A+p\wedge X,$$
and we get the decomposition 
$$\so (1, n + 1)_{\Real p} = (\mathR p\wedge q \oplus  \so(n))\zr p \wedge \mathR^n.$$ For the subalgebra of  $\so(1,n+1)$ annihilating the isotropic vector $p$ we have
$$\so (1, n + 1)_{p} =  \so(n)\zr p \wedge \mathR^n.$$

\section{Main Results}

Let $(M,g)$ be a Lorentzian manifold admitting a metric connection with parallel twistor-free torsion having non-vanishing vectorial component $\xi$. 

Assuming that $(M,g)$ is simply connected and complete and that $g(\xi,\xi)=\varepsilon=\pm 1$, one may apply all arguments from the proof of Theorem \ref{ThMP5.2} and conclude that $(M,g)$ must be isometric to the manifold given by \eqref{metricg}. But such  Lorentzian manifold is not complete \cite{ACGL,Bohle} and we see that  the statement of Theorem \ref{ThMP5.2} cannot be extended to complete Lorentzian manifolds with non-isotropic $\xi$. 
Since assuming geodesic completeness is usually too strong in Lorentzian geometry, one may consider weaker notions of geodesic completeness, such as space-like geodesic completeness time-like geodesic completeness, or light-like geodesic completeness (to be short, we use the notion of $\eta$-completeness, where $\eta=1,-1,0$). 
It turned out that assuming a weaker notion of geodesic completeness still does not allow to prove an analogue of Theorem \ref{ThMP5.2} for Lorentzian manifolds and finally we prove

\begin{theorem} \label{ThComplg}
Let $(M,g)$ be a Lorentzian manifold admitting a metric connection with parallel twistor-free torsion having non-vanishing vectorial component $\xi$. Then $(M,g)$ is not  $\eta$-complete for $\eta=1,-1,0$.
	\end{theorem}

Theorem \ref{ThComplg} generalizes the following two results.   Lorentzian space forms do not admit homogeneous structures of vectorial type~\cite{GO97}.  A homogeneous Lorentzian space admitting a homogeneous structure of isotropic vectorial type is  a singular homogeneous plane wave~\cite{MA01,BOL}.

We prove an analogue of Theorem \ref{ThMP5.3} for Lorentzian manifolds if we just require the completeness of the vectorial component of the torsion in the case when the vectorial component is non-isotropic.

\begin{theorem}\label{Th_nonis}
	Let $(M,g)$ be a simply connected Lorentzian manifold. Let $\varepsilon=\pm 1$. Then 
	 $(M,g)$  admits a metric connection with parallel twistor-free torsion having the  vectorial component given by a non-isotropic complete 
	vector field  
	 $\xi$ with $\varepsilon g(\xi,\xi)>0$ if and only if  $(M,g)$ is homothetic to the warped product \begin{equation}\label{metricg}(\Real\times N,\varepsilon (dt)^2+e^{2\varepsilon t}g_N),\end{equation}
	 where  $(N,g_N)$ is a Riemannian (if $\varepsilon=-1$) or a Lorentzian
	 (if $\varepsilon=1$) manifold carrying a parallel 3-form $\tau \in \Omega^3(N)$ that satisfies $\sigma_\tau= 0$.
	 \end{theorem}
	
	\begin{ex}
The flat slicing coordinates on an open subset of the de~Sitter space $dS^{n+2}$ allow one to write the de~Sitter metric in the form
$$g=-(dt)^2+e^{-2t}\sum_{i=1}^{n+1}(dy^i)^2.$$	
The vector field $\xi=\partial_t$ defines a metric connection $\nabla$ with vectorial torsion. This connection is a homogeneous structure on an open subset of $dS^{n+2}$ and, according to \cite{GO97}, $\nabla$ cannot be extended to the  entire $dS^{n+2}$. 
\end{ex}

Then we prove an analogue of Theorem \ref{ThMP5.2} for 
 the Lorentzian signature.

\begin{theorem}\label{ThLordT=0}
	Let $(N,g)$ be a    Lorentzian
	manifold carrying a parallel 3-form $\tau$ with $\sigma_\tau= 0$. Then $(N,g)$ is locally isometric to one of the following manifolds:
	\begin{itemize}
		\item[1.] $ \dim N = 3$, $(N,g)$ is an oriented 3-dimensional Lorentzian manifold and $\tau$ is a constant multiple of its Lorentzian volume form;
		\item[2.] $\dim N\geq 4$,   $(N,g)$ is a Lorentzian manifold carrying a parallel isotropic vector field $p$, and $\tau=p^\flat \wedge \omega$, where $\omega$ is a parallel 2-form on the screen bundle $p^\bot/\left<p\right>$.
		\item[3.] $(N,g)$ is a Lorentzian manifold and $\tau = 0$;
		\item[4.] $(N,g)$ locally is a product of a Lorentzian manifold
		 from one of the cases 1--3 with a Riemannian manifold from Theorem~\ref{ThMP5.2}.
	\end{itemize}
\end{theorem}

From Lemma \ref{lem} it follows that Theorem \ref{ThLordT=0} describes Lorentzian manifolds admitting metric connections $\nabla$ with closed $\nabla$-parallel skew-symmetric torsion.

Then we consider the case when the vectorial part of the torsion is isotropic. In that case the $g$-norm of the torsion coincides with the $g$-norm $||S||_g$  of the skew-symmetric component $S$. We say that the component $S$ is degenerate if $||S||_g=0$, otherwise we say that $S$ is non-degenerate. The first theorem deals with the case of degenerate skew-symmetric component $S$.

   \begin{theorem}\label{Th_isotropic_deg}
      Let $(M_0,g_0)$ be a  Lorentzian manifold  with a    $\nabla^{g_0}$-parallel isotropic vector field $p_0$ and a $\nabla^{g_0}$-parallel  form 
   $\tau_0=p_0^\flat\wedge \omega$, where $\omega$ is a  $\nabla^{g_0}$-parallel 2-form
   on the screen bundle $p^\bot/\left<p\right>$. 
      Suppose that there exists a function $\varphi$ on $M_0$ such that $p_0^\flat=de^{\varphi}$. Let $\xi=e^{-\varphi}p_0$ and $S=e^{-\varphi}\tau_0$. Then the connection $\nabla$ on $M_0$ given by \eqref{connectionnabla} has parallel 
      twistor-free torsion with isotropic vectorial component $\xi$ and the skew-symmetric component $S$.
      
      Conversely, each simply connected Lorentzian manifold   admitting a metric connection with parallel twistor-free torsion having isotropic vectorial component and degenerate skew-symmetric component is globally equivalent to a just constructed one.
         \end{theorem}

\begin{rem} Let $(M_0,g_0)$ be a Lorentzian manifold with a $\nabla^{g_0}$-parallel isotropic vector field $p_0$. 
	Then $d p^\flat_0=0$, and there exists a function $f$ such that $p^\flat_0=d f$. Such a function is defined up to a constant. The condition in the above construction is satisfied whenever each such function is bounded from below. For example, let $(B,b)$ be a Riemannian manifold. Consider the Walker manifold
	$$(M_0=\Real \times B\times \Real_+,\quad g_0=2dvdu+b+H(du)^2),$$
	where $v$ and $u$ are the coordinates on   $\Real$ and  $\Real_+$, respectively, and $H$ is any function on $B\times \Real_+$. Then  
	$p_0=\partial_v$ is a parallel isotropic vector, and $p^\flat_0 =du=de^{\varphi}$, $\varphi=\ln u$. The condition on $p^\flat_0$ from the above construction is satisfied.
\end{rem}

\begin{ex}
	Consider the metric
	$$g_0=2dvdu+\sum_{i=1}^n(dx^i)^2+(H+2v)(du)^2,\quad \partial_v H=0.$$
	The vector field $p_0=-e^{-u}\partial_v$ is $\nabla^{g_0}$-parallel and it satisfies $p_0^\flat=de^{-u}$.
	Let $$\xi=-\partial_v,\quad S=du\wedge \sum_{i<j}F_{ij}dx^i\wedge dx^j,$$
	where $F=(F_{ij})$ is a constant skew-symmetric matrix. Then the connection $\nabla$ given by \eqref{connectionnabla} is a metric connection with parallel twistor-free torsion having isotropic vectorial component $\xi$ and  skew-symmetric component $S$. The connection $\nabla$ is a homogeneous structure if and only if $\nabla R=0$. According to \cite{BOL,Meessen,MA01}, this is the case if and only if 
	$$H=A(e^{-uF}x,e^{-uF}x),$$ i.e, $g$ is a singular homogeneous plane-wave metric. Here $A$ is a constant symmetric bilinear form. 
\end{ex}

The following two theorems give a complete solution for non-degenerate $S$.

\begin{theorem}\label{Th_isotropic_nondeg_dim=3}
Let $(M,g)$ be a 3-dimensional Lorentzian manifold. Then 
$(M,g)$ admits a metric connection with parallel twistor-free torsion having isotropic vectorial component $p$ and non-degenerate skew-symmetric component $S$ if and only if $M$ is oriented and admits an isotropic  vector field $p$ 
such that \begin{equation}\label{condpKundt}\nabla^g_Xp=-g(X,p)p-\frac{1}{2}S(X,p),\quad\forall\,X\in\Gamma(TM),\end{equation}
where $$S=a\Vol_g,\quad a\in\Real,\quad a\neq 0.$$
In particular, such $(M,g)$ is a Kundt space, and locally $g$ may be written as
$$g=2dvdu+2avdxdu+(dx)^2+(-2ve^{-ax}+C(x,u))(du)^2,$$
where $C(x,u)$ is an arbitrary function and $p=e^{-ax}\partial_v$.
 \end{theorem}

\begin{theorem}\label{Th_isotropic_nondeg} Let $(M_0,g_0,p_0,\tau_0,\varphi)$ be a Lorentzian manifold as in
	 Theorem \ref{Th_isotropic_deg}.
 Let $(N,g_N)$ be a  Riemannian manifold with a non-zero $\nabla^{g_N}$-parallel 3-form $\tau_N$ satisfying $\sigma_{\tau_N}=0$. Consider the   manifold
$$M=M_0\times N$$
with the Lorentzian metric $$g=g_0+e^{2\varphi}g_N.$$
Let $$\xi=e^{-\varphi}p_0,\quad S=e^{-\varphi}\tau_0+e^{3\varphi}\tau_N.$$ Then the connection $\nabla$ on $M$ given by \eqref{connectionnabla} has parallel 
	twistor-free torsion with isotropic vectorial component $\xi$ and  non-degenerate skew-symmetric component $S$.
	
	Conversely, each Lorentzian manifold of dimension different from 3  admitting a metric connection with parallel twistor-free torsion having isotropic vectorial component and non-degenerate skew-symmetric component is locally equivalent to a just constructed one.
\end{theorem}

\section{Proof of Theorem \ref{Th_nonis}}

Let $(M,g)$ be a simply connected Lorentzian manifold. Suppose that 
$(M,g)$  admits a metric connection $\nabla$ with parallel twistor-free torsion having  vectorial component given by a non-isotropic complete 
vector field $\xi$. Let us assume that $g(\xi,\xi)=\varepsilon=\pm 1$. The connection $\nabla$ is given by \eqref{connectionnabla}. For the proof we use the ideas from \cite{MP}, and we apply \cite[Prop. 8]{LSch} to obtain the global decomposition of $M$.

		Consider the 1-form $\eta$ given by 
		$$ \eta(X) = \varepsilon g(\xi,X). $$
		It is clear that $\eta(\xi) = 1$. As in \cite{MP} it can be shown that 
		$$ d\eta = S(\xi) = 0.$$
		Proposition 8 from \cite{LSch} implies that all leaves of the foliation on $M$ tangent to the distribution $\ker\eta$ are pairwise diffeomorphic, and the manifold $M$ is diffeomorphic to the product
		$$M\cong\mathR \times N,$$
		where $N$ is a leaf of the foliation. Under this diffeomorphism, the vector field $\xi$ corresponds to $\partial_t$. From the definition of $\eta$ it follows that the vectors tangent to the leaves of the foliation are orthogonal to $\partial_t$. Consequently, the metric $g$ may be written in the form
		$$g=\varepsilon(dt)^2+h,$$
		where $h$ is a $t$-family of metrics on $N$.
		Consider the metric $$\tilde g=e^{-2\varepsilon t}g.$$
		Applying the formula for  the Levi-Civita connection under the conformal change, \eqref{connectionnabla}, and the fact that $S(\xi)=0$, we obtain
		$$\nabla^{\tilde g}_X\xi=-\varepsilon(Xt)\xi.$$
		This shows that the distribution generated by $\xi$ and the distribution $\ker\eta$ are  $\nabla^{\tilde g}$-parallel. The local version of the Wu Theorem implies that $$g_N=e^{-2\varepsilon t}h$$ is a metric on $N$ independent of $t$.
		Thus,
		$$g=\varepsilon (dt)^2+e^{2\varepsilon t}g_N.$$
		The equality $S(\xi)=0$ shows that $S$ is a $t$-family of 3-forms on the manifold $N$. As in \cite{MP} it can be shown that the condition $\nabla S=0$ is equivalent to the conditions
		$$S=e^{3\varepsilon t}\tau_N,\quad \tau_N\in\Omega^3(N),\quad \nabla^g\tau_N=0,\quad\sigma_{\tau_N}=0.$$
		This completes the proof of the theorem.
	\qed

\section{Proof of Theorem \ref{ThLordT=0}}
	Let $(N,g)$ be a   Lorentzian
	manifold carrying a non-zero parallel 3-form $\tau$ such that $\sigma_\tau= 0$. Consider the  metric connection
	$$\nabla=\nabla^g+\frac{1}{2}\tau$$
	with skew-symmetric torsion $\tau$.
	From Lemma \ref{lem} it follows that $\tau$ is $\nabla$-parallel. 
	In \cite{EG22} 
	we described holonomy, curvature and torsion of Lorentzian connections with parallel skew-symmetric torsion. We have now to consider the additional condition   $\sigma_\tau= 0$. We assume that $\dim N=n+2\geq 3$ and denote by  
	$\g\subset\so(1,n+1)$ the holonomy algebra of the connection $\nabla$.
	We fix a point $x\in N$. The tangent space $T_xN$ may be identified with the Minkowski space  $\Real^{1,n+1}$.

Recall that	a subalgebra $\g\subset\so(1,n+1)$ is called weakly irreducible if it does not preserve any proper non-degenerate subspace of $\Real^{1,n+1}$.  
	The geometry $(N, g, \nabla)$ is called \emph{reducible} if the holonomy algebra   $\g\subset\so(1,n+1)$ of the connection $\nabla$ is \emph{not} weakly irreducible, i.e., $\g$ preserves a proper non-degenerate subspace of the tangent space. In this case there exists a non-trivial $\g$-invariant orthogonal decomposition of the tangent space
	\begin{equation}\label{dec0}T_xM=L\oplus E.\end{equation}
	The geometry $(N, g, \nabla)$ is called \emph{decomposable}
	if the holonomy algebra   $\g\subset\so(1,n+1)$ preserves an orthogonal decomposition \eqref{dec0} such that it holds 
	$$\tau_x\in\wedge^3 L\oplus \wedge^3 E.$$ Otherwise we say that the geometry is \emph{indecomposable}. If the geometry is decomposable, then  it is a product of two other geometries, i.e., $(N,g)$ is as in the case 4 from the statement of the theorem. Thus we may assume that the geometry is indecomposable.
	
In Section 3 from \cite{EG22} we proved that if $\nabla$ is a metric connection on a Lorentzian manifold  $(N, g)$ with parallel skew-symmetric torsion $\tau$ and weakly irreducible holonomy algebra, then $\tau$ automatically satisfies the condition  
	$\sigma_\tau=0$, moreover, $(N,g)$ is as in the case 1 or 2 from the statement of the theorem.
	
Now we assume that the geometry $(N, g, \nabla)$ is reducible and indecomposable. Then the holonomy algebra $\g$ preserves a decomposition \eqref{dec0}. We may assume that the induced representation of $\g$ in $L$ is weakly irreducible. We consider several cases depending on the dimension of $L$ and use the description of $\tau$ from \cite{EG22}.
We will denote by $\tau$ also the value of the field $\tau$ at the point $x$.
Since the field $\tau$ is $\nabla$-parallel, it is enough to check that its value at the point $x$ satisfies the condition 	
		\begin{equation} \label{eq:tautau=0}
	\tau(X)\cdot \tau = 0 \quad \text{for all} \quad X \in T_xN.
	\end{equation}	If $\dim L\geq 2$, we denote by $p,e_1,\dots,e_k,q$ a Witt basis in $L$.	
	
	Let $\dim L = 1$. Then 
	\[ T_xM = \mathR e_- \oplus E, \]
	\[ \tau = e_- \wedge \theta + \omega_E, \] 
	where $g(e_-, e_-) = -1$,  $\theta \in \wedge^2 E$ and $\omega_E \in \wedge^3 E$ are annihilated by $\g$ and it holds $\theta \cdot \omega_E = 0$. The last condition  may be written in the form 
	\[\underset{X Y Z}\Su \omega_E(\theta(X),Y,Z) = 0 \quad \text{for all} \quad X,Y,Z \in E. \]
		Let $X \in E$. Then 
	\[ \tau(X) = \theta(X) \wedge e_- + \omega_E(X), \]
	\[ \tau(X)\cdot \tau = \theta(X) \wedge \theta + \omega_E(X) \cdot \omega_E + e_- \wedge (\omega_E(X)\cdot \theta) = 0. \]
	This is equivalent to 
	\begin{align}
	\theta(X) \wedge \theta + \omega_E(X) \cdot \omega_E &= 0, \label{eq:L1.1} \\
	\omega_E(X)\cdot \theta &= 0. \label{eq:L1.2}
	\end{align}
	We have 
	\begin{multline*}
	0 = \omega_E(X)\cdot \theta(Y,Z) = -\theta(\omega_E(X,Y),Z) - \theta(Y,\omega_E(X,Z)) = \\ 
	= \omega_E(X,Y,\theta(Z)) - \omega_E(X,Z,\theta(Y)) = \omega_E(\theta(Z),X,Y) + \omega_E(\theta(Y),Z,X) = \\ 
	= - \omega_E(\theta(X),Y,Z),
	\end{multline*}
	which means that $\theta \in \wedge^2(\ker \omega_E)$.  Since the geometry is indecomposable, this implies that $\omega_E = 0$.
	Equation \eqref{eq:L1.1} now reads as
	\[ \theta(X) \wedge \theta = 0, \]
	 for all $X\in E$. This is possible only if $\operatorname{rk} \theta \leq 2$. 
	Thus the indecomposability implies that the dimension of $N$ is 3.
	
	Let $\dim L = 2$. In this case $\tau$ has the form 
	\[ \tau = p \wedge q \wedge v + \theta \wedge v + \omega_{E_1}, \]
	where $v\in E$ is a non-zero vector, $E = \mathR v \oplus E_1$ is an orthogonal decomposition, $\g v = 0$, $\omega_{E_1} \in \wedge^3 E_1$, $\theta \in \wedge^2 E_1$, $\g\cdot\omega_{E_1}=0$, $\g\cdot\theta=0$. 
	Consider the condition \eqref{eq:tautau=0}. 
	It holds
	\[ \tau(p) = -p \wedge v, \]
	\[ \tau(p)\cdot \tau = g(v,v)\theta \wedge v = 0. \]
Hence, $\theta = 0$.  The indecomposability again implies that the dimension of $N$ is 3.
	
	Assume that $\dim L = 3$. Then, given an arbitrary Witt basis $p, e_1, q$ of $L$, the torsion $\tau$ has the form
	\[ \tau = p \wedge (\alpha e_1 \wedge q + e_1 \wedge v + \lambda) + \theta \wedge v + \omega_{E_1},\]
	where $v\in E$ is a  vector, $E_1$ is the orthogonal complement of $v$ in $E$,  $\g v = 0$, $\omega_{E_1} \in \wedge^3 E_1$, $\theta \in \wedge^2 E_1$ and $\lambda \in \wedge^2 E$ are annihilated by $\g$.
	It holds that
	\[ \tau(e_1) = -\alpha p \wedge q - p \wedge v,\]
	\[ \tau(e_1) \cdot \tau = p \wedge (\alpha \lambda + g(v,v) \theta) = 0.\]
	Consequently, 
	\[ \alpha \lambda + g(v,v) \theta = 0. \]
	
	First suppose that $\alpha\ne 0$. Let us set 
	\[
	\tilde v = v - \frac{g(v,v)}{\alpha} p, \quad 
	\tilde q = q + \frac{1}{\alpha}\left(v - \frac{1}{2\alpha}g(v,v)p\right), \quad
	\tilde E = \mathR \tilde v \oplus E,   \quad 
	\tilde L = \left< p, e_1, \tilde q \right>.
	\]
	Then $p,e_1,\tilde q$ is a Witt basis in $\tilde L$ and $T_xM = \tilde L \oplus \tilde E$ is a new $\g$-invariant orthogonal decomposition. 
	The torsion $\tau$  now may be written as
	\[ \tau = \alpha p \wedge e_1 \wedge \tilde q + \theta \wedge \tilde v + \omega_{E_1} \]
	and we see that the geometry is decomposable.
	
	Suppose now that $\alpha = 0$. The condition \eqref{eq:tautau=0} may be rewritten in the following way:
	\[
	\theta\wedge v = 0, \quad
	\lambda(v) = 0, \quad
	\lambda \cdot \omega_{E_1} = 0, \quad
	\omega_{E_1}(X) \cdot \omega_{E_1} = 0,\quad\forall X\in E_1.
	\]
	The expressions for the curvature tensor from \cite[Theorem 5]{EG22} imply that 
	\[ \g = \mathR p \wedge e_1 \oplus \mathfrak{b}, \]
	where $\mathfrak{b}$ is the projection of $\g$ onto $\so(E_1)$. 	
	The above equalities show that $\omega_{E_1}$ defines a Lie algebra structure on $E_1$,  and $\lambda$ is a derivation of the Lie algebra $E_1$. 
	Let $E_0=\ker \omega_{E_1}\subset E_1$ and let $E_1'$ be the orthogonal complement to $E_0$ in $E_1$. Then $E_0$ is a commutative Lie algebra, and $E_1'$ is either the trivial or a semisimple Lie algebra. 
		Therefore,
	\[ \lambda = \omega_{E_1}(U_0) + \lambda_0, \]
	where $U_0 \in E'_1$ and $\lambda_0 \in \wedge^2 E_0$.
	Since $\g$ annihilates $\lambda$ and $\omega_{E_1}$, it annihilates the vector $U_0$.
	
	Let 
	\[
		\tilde E = \{ \tilde U = U + g(U,U_0)p \, | \, U \in E \},\quad \tilde q = q - U_0 - \frac{1}{2}g(U_0,U_0)p,\quad
	\tilde L = \left< p, e_1, \tilde q \right>.
	\] 
	Then,  $T_xM = \tilde L \oplus \tilde E$ is a $\g$-invariant orthogonal decomposition. Now,
	$$\tau=p \wedge (e_1 \wedge \tilde v + \tilde\lambda_0) +  \omega_{\tilde E_1},$$
	where in the new notation $\tilde E=\Real \tilde v+\tilde E_1$, $\tilde E_1=\tilde E_0\oplus \tilde E'_1$, $\tilde\lambda_0\in\wedge^2 \tilde E_0$, $\omega_{ \tilde E_1}\in\wedge^3 \tilde E'_1$. The indecomposability implies that $\tilde E'_1=0$. Thus, $(N,g)$ is as  in the case 1 or 2 from the statement of the theorem.

Finally let us assume that  $\dim L \geq 4$. Then the  torsion $\tau$ has the form
	\begin{equation}
	\label{eq:tau}
	\tau = p \wedge (\omega + \sum_{i=1}^k e_i \wedge \mu_i + \lambda) + \omega_E,
	\end{equation}
	where $\mu_i \in E$, $\omega_E \in \wedge^3 E$, $\omega \in \wedge^2 \mathR^k$, $\lambda \in \wedge^2 E$, and it holds that $\lambda \cdot \omega_E = 0$.
	Let us consider the condition \eqref{eq:tautau=0}:
	\[ \tau(e_i) = -p \wedge \mu_i - p \wedge \omega(e_i), \]
	\[ \tau(e_i) \cdot \tau = -p \wedge \omega_E(\mu_i) = 0, \]
	\[ \tau(U)\cdot \tau = \omega_E(U) \cdot \omega_E = 0,\quad U\in E.
	 \]
The condition $\tau(q) \cdot \tau = 0$ is equivalent to $\lambda \cdot \omega_E = 0$. 
	From these equalities we get $\omega_E(\mu_i) = 0$. 
	Again, $\omega_E$ defines a  Lie algebra structure on $E$, and as in the previous case, the indecomposability of the geometry implies that $\omega_E=0$, and  $(N,g)$ is as  in the case 1 or 2 from the statement of the theorem. 
	\qed

\section{Proof of Theorem \ref{Th_isotropic_nondeg_dim=3}}

	Let $(M,g)$ be a 3-dimensional  Lorentzian manifold. 
Suppose that $(M,g)$ admits a metric connection $\nabla$ with parallel twistor-free torsion having isotropic vectorial component $\xi$ and non-zero skew-symmetric component $S$. 

The connection $\nabla$ is given by \eqref{connectionnabla}.
Since $S$ is parallel, it is non-vanishing, and, consequently, $M$ is orientable. Fix an orientation of $M$, and let $\Vol_g$ be the volume form defined by $g$. The torsion is $\nabla$-parallel if and only if $\nabla p=0$ and $\nabla S=0$.

It holds $$0=\nabla S=\nabla^g S,$$ i.e., $S$ is proportional to the volume form, 
$$S=a \Vol_g,\quad a\in\Real,\quad a\neq 0.$$
Now, the torsion is parallel
if and only if
$\nabla p=0$. By \eqref{connectionnabla}, the last condition  is equivalent to  \eqref{condpKundt}.
It is obvious that the vector field $p$ satisfies the conditions \eqref{Kundtfield}, i.e., $(M,g)$ is a Kundt space.
Consequently, the metric $g$ may be locally written in the form
$$g=2dvdu+2A(v,x,u)dxdu+(dx)^2+H(v,x,u)(du)^2.$$
The vector field $p$ satisfies 
$p=\varphi\partial_v$, for a function $\varphi=\varphi(x,u)$.
Consider  the local field of Witt frames 
$$p=\varphi\partial_v,\quad e=\partial_x,\quad q=\frac{1}{\varphi} \left(\partial_u-A\partial_x-\frac{1}{2}H\partial_v\right).$$
The condition \eqref{condpKundt} is equivalent to
$$\nabla^g_p p=0,\quad \nabla^g_e p=-\frac{1}{2}ap,\quad \nabla_q^g p=-p  +\frac{1}{2}ae.$$
It is easy to check that
$$\nabla^g_e p=\left(\partial_x\varphi+\frac{1}{2}\varphi\partial_vA\right)\partial_v.$$
This implies that
$$\frac{\partial_x\varphi}{\varphi}+\frac{1}{2}\partial_vA=-\frac{1}{2}a.$$
Next, $$\nabla^g_q p=\left(\frac{1}{\varphi}(\partial_u\varphi-A\partial_x\varphi)+\frac{1}{2}\partial_vH-A\partial_vA\right)\partial_v+\frac{1}{2}\partial_vA\partial_x.$$
We conclude that $\partial_v A=a$, i.e., $A=av+B(x,u)$, and $\partial_x\varphi=-\varphi,$ i.e., $\varphi=c(u)e^{-ax}$. A simple coordinate transformation allows to assume that $B(x,u)=0$ and $c(u)=1$.
Finally we get that $H=-2e^{-ax}v+C(x,u)$.  \qed

\section{Proof of Theorems \ref{Th_isotropic_deg} and \ref{Th_isotropic_nondeg}}\label{SecThAB}

Let $(M,g)$ be a Lorentzian manifold admitting a metric connection $\nabla$ with parallel twistor-free torsion having isotropic vectorial component $p$ and skew-symmetric component $S$. If $\dim M=3$, then thanks to Theorem \ref{Th_isotropic_nondeg_dim=3} {\it we will assume that} $||S||_g=0$; since in dimension 3 each 3-form is proportional to the (local) volume form, this implies that $S=0$. 
Denote by $\g\subset \so(1,n+1)$ the holonomy algebra of the connection $\nabla$ at a point $x\in M$. By  abuse of notation, we denote by $p$, $S$, $R$ the values of the corresponding tensor fields at the point $x$.
Then $\g$ preserves the vector $p\in\Real^{1,n+1}=T_xM$, i.e., $\g\subset\so(1,n+1)_{p}$.

\begin{lem}\cite[Lemma 1]{EG22}
If $\g\subset \so(1,n+1)_p$ is a weakly irreducible subalgebra, then each $S\in\wedge^3\Real^{1,n+1}$ annihilated by $\g$ is of the form  $S=p\wedge \omega$
for a bivector $\omega$ on $p^\bot/\left<p\right>$.
\end{lem}

\begin{lem}\label{LemnondegS} If the holonomy algebra $\g\subset \so(1,n+1)$ of the connection $\nabla$ is not weakly irreducible, then
	$\g$ preserves an orthogonal decomposition $$\Real^{1,n+1}=L\oplus E,$$ where $$ L=\Real p\oplus\Real^k\oplus\Real q,\quad  0\leq k\leq n,$$
	such that 
	$$S=p\wedge \omega+\omega_E,$$
	where $\omega\in\wedge^2\Real^k$ and $\omega_E\in\wedge^3 E$.
	Moreover, $\sigma_{\omega_E}=0$ and $\ker \omega_E=0$.
\end{lem}

Note that the statement of the lemma includes the case $E=0$.

{\bf Proof of Lemma \ref{LemnondegS}.} 
If $\g\subset \so(1,n+1)$ is not weakly irreducible, then it is clear that
$\g$ preserves an orthogonal decomposition $$\Real^{1,n+1}=L\oplus E,\quad \dim L\geq 1$$
such that the induced representation of $\g$ in $L$ is weakly irreducible.

First assume that $\dim L\geq 4$.
By \cite[Lemma 7]{EG22}, $$ S = p \wedge \omega + \omega_E,$$ where 
$$ \omega = \omega_{\Real^k} + \mu + \lambda,$$
here $\omega_{\Real^k} \in \wedge^2 \mathR^k$, $\mu \in \mathR^k \wedge E$, $\lambda \in \wedge^2 E$.
The Bianchi identity \eqref{Bidcomp} written for the vectors $U,V\in E$ and $q$ takes the form
$$R(U,V)q+R(V,q)U+R(q,U)V=-(\omega\cdot\omega_E)(U,V)+\lambda(U,V)p+\omega_E(U,V).$$
Since $E$ is holonomy-invariant, $R(V,q)U,R(q,U)V\in E$.
Hence, multiplying the above equality scalarly by $q$, we get
$$\lambda=0.$$ This implies that $\omega\cdot\omega_E=0$, and
$$R(V,q)U+R(q,U)V=\omega_E(U,V).$$
Since the first prolongation of $\so(E)$ is trivial, the equation 
$R(V,q)U+R(q,U)V=0$, for all $U,V\in E$, has only trivial solution. This shows that 
$$R(q,U)V=\frac{1}{2}\omega_E(U,V),\quad\forall\, U,V\in E.$$
Since $\g$ annihilates $S$, and $R$ takes values in $\g$, we get
$$(R(q,U)\cdot S)(U_1,U_2,U_3)=0,\quad \forall\, U,U_1,U_2,U_3\in E.$$
This implies that
$$\omega_E(U)\cdot \omega_E=0,\quad \forall\, U\in E,$$
i.e.,
$$\sigma_{\omega_E}=0.$$
We conclude that the equality
$$g([U_1,U_2],U_3)=\omega_{E}(U_1,U_2,U_3),\quad \forall\, U_1,U_2,U_3\in E$$ defines a Lie bracket on the vector space $E$.
Let $E_0=\ker \omega_E$, and let $E'$ be the orthogonal complement to $E_0$ in $E$. Then $E_0$ and $E'$ are  commutative and semisimple ideals in $E$, respectively. 

The Bianchi identity written for the vectors  $X\in \Real^k$, $U\in E$ and $q$ easily implies that $g(R(q,U)X,Y)=0$ for all $Y\in\Real^k$. Hence, $R(q,U)X\in\Real p$.  Let us write the bivector $\mu$ in the form $\mu=\sum_{i=1}^k e_i\wedge \mu_i$, where $\mu_i\in E$. Now, $\g$ annihilates the tensor $p\wedge\mu$, hence,
$$0=R(q,U)\cdot(p\wedge \mu)=\sum_{i=1}^k p\wedge e_i\wedge R(q,U)\cdot\mu=\frac{1}{2}\sum_{i=1}^k p\wedge e_i\wedge[U,\mu_i].$$
This implies that $[U,\mu_i]=0$ for all $U\in E$. Consequently,  $\mu_i\in E_0$. Denoting now $L\oplus E_0$ by $L$ and $E'$ by $E$, we see that $S$ is just as in the statement of the lemma.

Assume now again that the representation of $\g$ in $L$ is weakly irreducible.

 Suppose that $\dim L=3$. According to \cite{EG22}, $S$ has the form
\[ S = p \wedge \omega + \omega_E + \alpha p \wedge e_1 \wedge q, \]
where $\alpha \in \mathbb{R}$ and all other elements are as above. 
Using the Bianchi identity for the  vectors $U\in E$, $e_1$ and $q$, it is easy to shown that $\alpha = 0$, and further considerations are just as in the previous case. 

Since $\g$ annihilates the isotropic vector $p$, it is not possible that the representation  
 of $\g$ in $L$ is weakly irreducible and $\dim L=2$. Finally, if $\dim L=1$, then since $\g$ annihilates the isotropic vector $p$, it annihilates a two-dimensional Lorentzian subspace of $\Real^{1,n+1}$. Again, it is not hard to show that the statement of the lemma holds true. The lemma is proved. \qed
 
 \begin{cor} Under the current assumptions it holds $S(p)=0$ unless $\dim M=3$ and $(M,g)$ is as in Theorem \ref{Th_isotropic_nondeg_dim=3}.  \end{cor}
 
Let us suppose  that the manifold $M$ is simply connected.
 Consider the dual 1-form $\eta$ to the vector field $p$,
 $$\eta(X)=g(p,X),\quad\forall\, X\in\Gamma(TM).$$
 It holds 
 \begin{multline*} d\eta(X,Y) = X g(p,Y) - Y g(p,X) - g(p,[X,Y]) \\= 
 g(\nabla^g_Xp,Y)+g(p,\nabla^g_XY)-g(\nabla^g_Yp,X)-g(p,\nabla^g_YX)-
 g(p,\nabla^g_XY-\nabla^g_YX)\\= 
 g(\nabla^g_X p, Y) - g(\nabla^g_Y p, X) =g\left(-(X\wedge p)p-\frac{1}{2}S(X,p),Y\right)-g\left(-(Y\wedge p)p-\frac{1}{2}S(Y,p),X\right)= 0. \end{multline*}
  Then there exists a function $\varphi$ such that $\eta=d\varphi$.
   
  {\bf Proof of Theorem \ref{Th_isotropic_deg}.} Suppose now that $S$ is degenerate. From the above lemmata it follows that $S=p\wedge \omega$.
  
  From \eqref{connectionnabla} it follows that
  $$\nabla^g_Xp=-g(X,p)p=-\eta(X)p.$$
  Consequently the vector field $$p_0=e^{\varphi}p$$
  is $\nabla^g$-parallel. Its dual 1-form $p_0^\flat$ is closed, and there exists a function $u$ on $M$ such that $du=p_0^\flat$. We conclude that 
  $$d\varphi=e^{-\varphi}du,$$ i.e.,
  $$d e^{\varphi}=du.$$ Since the both functions $\varphi$ and $u$ are defined up to a constant, we may assume that
  $$u=e^\varphi,$$
  i.e., the function $u$ is positive. By \eqref{connectionnabla}, $\nabla^g \tau_0=0$, where $\tau_0=e^{\varphi} S$.  \qed

  {\bf Proof of Theorem \ref{Th_isotropic_nondeg}.}
Suppose that the manifold $M$ is simply connected. Suppose that $||S||_g\neq 0$ and $\dim M\neq 3$. Let $L$ and $E$ be as in Lemma \ref{LemnondegS}. Since the subspaces $L$ and $E$ of $T_xM$ are $\g$-invariant, they define  $\nabla$-parallel distributions $\mathcal{L}$ and $\mathcal{E}$ on $M$. Moreover, $$S=S_1+S_2,\quad S_1\in \wedge^3 \Gamma(\mathcal{L}),\quad S_2\in \wedge^3 \Gamma(\mathcal{E}),\quad \nabla S_1=\nabla S_2=0,\quad \sigma_{S_1}=\sigma_{S_2}=0.$$
  Consider the new metric   $$ h = e^{-2\varphi}g,$$ where the function $\varphi$ is as above.

  By the standard formula,  for all $X,Y \in \Gamma(TM)$ it holds
  \[ 
  \nabla^h_X Y = \nabla^g_X Y - g(p,X)Y - g(p,Y)X + g(X,Y)p. 
  \]
  Combining this equality with \eqref{connectionnabla}, we get 
  \begin{equation}\label{nabla_h}
  \nabla^h_X Y = \nabla_X Y - g(p,X)Y - \frac{1}{2}S(X,Y),
  \end{equation}
  or just 
   \begin{equation}\label{nabla_h'}
  \nabla^h_X  = \nabla_X  - g(p,X)\id - \frac{1}{2}S(X).
  \end{equation}
  This implies that the distributions $\mathcal{L}$ and $\mathcal{E}$
 are $\nabla^h$-parallel. By the Wu Theorem, $(M,h)$ is locally a product of a Lorentzian manifold $(M_0,h_0)$ and a Riemannian manifold $(N,g_N)$, where $M_0$ and $N$ are integral submanifolds of the distributions $\mathcal{L}$ and $\mathcal{E}$, respectively. In particular, locally it holds
 $$h=h_0+g_N.$$ 
 
 The vector field $p$ is tangent to the distribution $\mathcal{L}$, i.e., $p$ is a family of vector fields on the manifold $M_0$ depending on the local coordinates on the manifold $N$. If $U\in\Gamma(TN)$, then $\nabla^h_Up=0$. Thus, $p$ is a vector field on $M_0$. If $X\in\Gamma(TM_0)$, then $$\nabla^{h_0}_Xp=\nabla^h_Xp=-g(X,p)p=-\eta(X)p=-d\varphi(X) p.$$
 This shows that $\varphi$ is a function on the manifold $M_0$.
 Consequently, $$g_0=e^{2\varphi}h_0$$ is a Lorentzian metric on $M_0$ and it coincides with $g$ restricted to $M_0$. We get that
 $$g=g_0+e^{2\varphi}g_N.$$ If $U\in\Gamma(TN)$, then by \eqref{nabla_h'}, $$\nabla^h_US_1=\nabla_US_1=0,$$ i.e., $S_1$ is a 3-form on $M_0$.
 Next, if $X\in \Gamma(TM_0)$, then, by \eqref{connectionnabla}, $$\nabla_X^{g_0}S_1=\nabla_X^{g}S_1=\nabla_XS_1=0.$$
 We conclude that the data $$(M_0,g_0,p_0=e^{\varphi}p,\tau_0=e^{\varphi}S_1,\varphi)$$ is as in Theorem \ref{Th_isotropic_deg}.
 
 For $U\in\Gamma(TN)$ it holds
 $$\nabla^{g_N}_US_2=\nabla^h_U S_2=\nabla_U S_2=0.$$
 If $X\in\Gamma(TM_0)$, then
 $$\nabla^h_X S_2=-g(p,X)\id\cdot S_2=3g(p,X) S_2=3d\varphi(X)S_2.$$
 Consequently, $$\nabla^h_X\tau_N=0,$$
 where $\tau_N=e^{-3\varphi}S_2$. Thus,
 $$S=e^{-\varphi}\tau_0+e^{3\varphi}\tau_N,\quad \nabla^{g_N}\tau_N=0,\quad \sigma_{\tau_N}=0.$$
 This proves the theorem. \qed

\section{Proof of Theorem~\ref{ThComplg}}

Let $(M,g)$ be a Lorentzian manifold that admits a metric connection $\nabla$ with parallel twistor-free torsion having non-zero vectorial component $\xi$. 
Passing to the universal covering, we may assume that the manifold is simply connected.

First suppose that $g(\xi,\xi)=0$. Then $(M,g)$ is given by one of the Theorems \ref{Th_isotropic_deg},  \ref{Th_isotropic_nondeg_dim=3}, \ref{Th_isotropic_nondeg}.	

Let $(M,g)$ be as in Theorem \ref{Th_isotropic_nondeg_dim=3}. 
	Let $\gamma(t)$ be a geodesic such that \begin{equation}\label{geodgam} g(\dot\gamma(0),p)=a\neq 0 \text{ and } g(\dot\gamma(0),\dot\gamma(0))=\eta.\end{equation} It is clear that  $\eta$ may be chosen to take any of the values $1,-1,0$.
	Consider the function
	$\alpha(t)= g(\dot\gamma(t),p)$
	defined along the geodesic $\gamma(t)$.
	It holds
	\begin{multline*}\dot\alpha(t)=\frac{d}{dt}\alpha(t)=\dot\gamma(t)\alpha(t)=
	\dot\gamma(t)g(\dot\gamma(u),p)=g(\nabla^g_{\dot\gamma(t)}\dot\gamma(t),p)+g(\dot\gamma(t),\nabla^g_{\dot\gamma(t)}p)\\=
	g\left(\dot\gamma(t),-g(\dot\gamma(t),p)p-\frac{1}{2}S(\dot\gamma(t),p)\right)=-g(\dot\gamma(t),p)^2-\frac{1}{2}S(\dot\gamma(t),p,\dot\gamma(t))=-\alpha^2(t).\end{multline*}
	Thus,
	$$\dot\alpha(t)+\alpha^2(t)=0.$$
	This shows that $$\alpha(t)=\frac{1}{t+c},\quad c\in\Real$$
is not defined for all values $t\in \Real$, i.e.,  geodesic $\gamma$ is not complete.

	Let now  $(M,g)$ be as in Theorem \ref{Th_isotropic_deg} or   \ref{Th_isotropic_nondeg}, i.e., $(M,g)$ is as in Section \ref{SecThAB}. We use the notation of this section.
	If $\dim M=3$, then we may assume that $S=0$. 
	 Let $\gamma(t)$ be a geodesic satisfying~\eqref{geodgam} with $\eta$ having one of the values $1,-1,0$. Since $p_0$ is $\nabla^g$-parallel, for each $t\in\Real$ it holds $g_{\gamma(t)}(\dot\gamma(t),(p_0)_{\gamma(t)})=a$. Let $\alpha(t)=\varphi(\gamma(t))$ be a function along $\gamma$, where $\varphi$ is as in Section \ref{SecThAB}. It holds
			\[ \frac{d}{dt}\alpha(t) =\dot\gamma(t)\varphi= d\varphi(\dot\gamma(t)) = g(\dot\gamma(t), p) = e^{-\alpha(t)} g(\dot\gamma(t),  p_0) = ae^{-\alpha(t)}. \]
	We conclude that
	$$ e^{\alpha(t)} = at + c, \quad c \in \mathR, $$
	and wee see that $\alpha(t)$ cannot be defined for all $t \in \mathR$, i.e., $\gamma(t)$ is not complete.

	Next we suppose that $g(\xi,\xi)\neq 0$. We may assume that $g(\xi,\xi)=\varepsilon\pm 1$.
	As in the proof of Theorem \eqref{Th_nonis}, it can be shown that $S(\xi)=0$.
	Let as above, $\gamma(t)$ be a geodesic satisfying
	$$  g(\dot\gamma(0),\xi)\neq 0 \text{ and } g(\dot\gamma(0),\dot\gamma(0))=\eta $$ with $\eta$ having one of the values $1,-1,0$. As above, it is easy to show that  the function $\alpha(t)=g(\dot\gamma(t),\xi)$ defined along $\gamma(t)$ satisfies the equation 
	$$\dot\alpha(t)+\alpha^2(t)=\eta\varepsilon.$$
	As we have seen, if $\eta=0$, then the geodesic $\gamma(t)$ is not complete.
	If $\eta\varepsilon=-1$, then $\alpha(t)$ satisfies
	$$\arctan(\alpha(t))=-t+c.$$
	Since the function $\arctan$ is bounded, $\gamma(t)$ is not complete.
	
	Suppose that $\eta\varepsilon=1$, i.e., $\eta=\varepsilon$.
		As in the proof of Theorem \ref{Th_nonis}, it can be shown that
	  $\nabla^g_\xi\xi=0$, i.e., $\xi$ is a geodesic vector field. We conclude that $\xi$ is a complete vector field. This allows us to apply Theorem \ref{Th_nonis} and to conclude that $(M,g)$ is isomorphic to the warped product
	  $$(\Real\times N, \varepsilon(ds)^2+e^{2s}g_N),$$
  where  $(N,g_N)$ is a Riemannian (if $\varepsilon=-1$) or a Lorentzian
(if $\varepsilon=1$) manifold.
A curve $\gamma(t)$ in $M$ is defined by a function $s(t)$ and a curve $\gamma_1(t)$ in $N$. According to \cite{ACGL,Bohle}, if $\gamma(t)$ is a geodesic, then the function $s(t)$ satisfies the equation
	  $$\ddot s=\varepsilon e^{2s}g_N(\dot\gamma_1,\dot\gamma_1).$$
	  Next, $$g(\dot\gamma,\dot\gamma)=\varepsilon \dot s^2+e^{2s}g_N(\dot\gamma_1,\dot\gamma_1).$$
	  This implies 
	  $$  \ddot s+\dot s^2= \varepsilon g(\dot\gamma,\dot\gamma).$$ Let $a>0$.
	  Consider a geodesic $\gamma(t)$ such that $g(\dot\gamma(0),\dot\gamma(0))=\varepsilon a^2$ and $\dot s(0)>a$. 
	  We obtain the equation $$  \ddot s+\dot s^2= a^2.$$
	  If $\dot s(t_0)=a$ for some $t_0\in\Real$, then $s(t)=at+c$ is the unique solution satisfying the conditions $s(t_0)=at_0+c$, $\dot s(t_0)=a$. Then
	  $\dot s(0)=a$, which contradicts the assumption $\dot s(0)>a$. We conclude that $\dot s(t)>a$ for all $t$.    
	Then $$\frac{\ddot s}{\dot s^2-a^2}=-1,$$
	and $$\ln\left|\frac{\dot s+a}{\dot s-a}\right|=2a(t-c),\quad c\in\Real.$$
	Since $\dot s(t)>a$ for all $t$, the function $\dot s(t)$ is not defined for $t=c$, and the geodesic $\gamma(t)$ is not complete. 
		 	\qed
	
	Theorems \ref{ThComplg} and  \ref{Th_nonis} imply
	
	\begin{cor} The Lorentzian warped product \eqref{metricg} is not  $\eta$-complete for $\eta=1,-1,0$.  \end{cor}

\vskip0.5cm

{\bf Acknowledgements.} The authors are thankful to Andrei Moroianu, Thomas Leistner and Eivind Schneider for useful discussions. The authors are grateful to the Reviewer for the important comments. I.E. was supported by  grant  MUNI/A/1099/2022 of Masaryk University.  A.G. acknowledges institutional support of University of Hradec Kr\'alov\'e.


\begin{thebibliography}{10}
	
	
	\bibitem{AF04}  I.\,Agricola, Th.\,Friedrich, On the holonomy of connections with skew-symmetric torsion. Math. Annalen 328 (2004), no. 4, 711-748.
	
	
		\bibitem{A10}  I.\,Agricola,
	Non-integrable geometries, torsion, and holonomy. in Handbook of pseudo-Riemannian Geometry and Supersymmetry,
	IRMA, EMS, 2010, 277--346.
	
\bibitem{AFeF15} I. Agricola, C. Ferreira, Th. Friedrich, 	The classification of naturally reductive homogeneous spaces in dimensions $n\leq 6$. 	Diff. Geom. Appl. 39 (2015), 59--92.

\bibitem{ACGL} D. Alekseevsky, V. Cortes, A. Galaev, Th.  Leistner, Cones over pseudo-Riemannian manifolds and their holonomy. Journal fur die Reine und Angewandte Mathematik 635 (2009), 23-69. 

\bibitem{BOL} 	M. Blau, M. O'Loughlin, Homogeneous plane waves. Nuclear Phys. B, 654 (2003), no. 1-2, 135--176.

	
\bibitem{Bohle}	C. Bohle, Killing spinors on Lorentzian manifolds. J. Geom. Phys. 45 (2003), 285–-308.
	
	
\bibitem{Walkerbook} M. Brozos-V\'azquez, E. Garc\'ia-R\'io, P. Gilkey, S. Nik\v{c}evi\'c, R. V\'azquez-Lorenzo, The
geometry of Walker manifolds, Synth. Lect. Math. Stat., 5, Morgan \&
Claypool
Publishers, Williston, VT, 2009.

\bibitem{CLBook} G. Calvaruso, M.C. L\'opez, Pseudo-Riemannian Homogeneous Structures. Springer 2019.

	\bibitem{CMS21}
R. Cleyton, A. Moroianu, U. Semmelmann,
Metric connections with parallel skew-symmetric torsion. Adv. Math. 378 (2021) 107519.

\bibitem{CHPP}A. Coley, S. Hervik, G. Papadopoulos, and N. Pelavas,
 Kundt spacetimes. Class. Quantum Grav. 26 (2009), no. 10, arc. num. 105016.

\bibitem{EG22}   I.\,Ernst, A.S.\,Galaev, On Lorentzian connections with parallel skew torsion. Documenta Mathematica 27 (2022), 2333--2383.


 
 	\bibitem{FFSPPM} {J. Figueroa-O'Farrill,  S. Philip, and P. Meessen}, Homogeneity and plane-wave limits. J. High Energy Phys. 05(05) (2005). 
 
 
 \bibitem{FI02}  Th. Friedrich, S. Ivanov, Parallel spinors and connections with skew-symmetric torsion in string theory. Asian J. Math. 6 (2002), no. 2, 303--335.


\bibitem{GO97} P.M. Gadea, J.A. Oubi\~na, Reductive homogeneous pseudo-Riemannian manifolds. Monatsh. Math. 124 (1997), 17--34.

	\bibitem{LScr} Th.~Leistner, Screen bundles of Lorentzian manifolds and some generalisations of pp-waves. J. Geom. and Phys. 56 (2006), no. 10, 2117--2134.


\bibitem{LSch}	T.\,Leistner, D.\,Schliebner, Completeness of compact Lorentzian manifolds with abelian holonomy. Math. Annalen  364 (2016), 1469--1503. 


\bibitem{Meessen} P. Meessen, 
Homogeneous Lorentzian spaces admitting a
homogeneous structure of type $T_1 \oplus T_3$. J. Geom. Phys. 56 (2006) 754--761.

\bibitem{MSh} \'A.\,Murcia, C.S.\,Shahbazi,
Contact metric three manifolds and Lorentzian geometry with torsion in six-dimensional supergravity.
J.  Geom. Phys. 158 (2020),
103868.

\bibitem{MA01} A. Montesinos Amilibia, Degenerate homogeneous structures of type $S_1$ on pseudo-Riemannian manifolds. Rocky Mountain J. Math. 31 (2001), 561--579.	


\bibitem{MP}	A.\,Moroianu, M.\,Pilca,	Metric connections with parallel twistor-free torsion. Int. J. Math. 32 (2021), arc. num. 2140011.

	
	

\bibitem{Str}  A. Strominger, Superstrings with torsion. Nucl. Phys. B274 (1986), 253--284.

		
\bibitem{TV}	F. Tricerri and L. Vanhecke. Homogeneous structures on Riemannian manifolds, volume 83 of
	London Mathematical Society Lecture Note Series. Cambridge University Press, Cambridge,
	1983.
	
	
	


\end{thebibliography}
\end{document}